\documentclass[12pt]{article}
\usepackage{amssymb,latexsym,amsmath}
\def\C{\mathbb C}
\def\R{\mathbb R}

\newtheorem{thm}{Theorem}[section]
\newtheorem{lem}{Lemma}[section]

\newtheorem{prop}{Proposition}[section]

\begin{document}

\sffamily
\title{Non-real zeros of  derivatives in the unit disc}
\author{J.K. Langley}
\maketitle

\begin{abstract}
The main result establishes an estimate for the growth of 
a real meromorphic
function $f$ on the unit disc $\Delta$  such that: (i)
 at least one of  $f$ and $1/f$ has finitely many poles and non-real zeros in $\Delta$;
(ii)~$f^{(k)}$ has   finitely many non-real zeros in $\Delta$, for some $k \geq 2$. \\
Keywords: meromorphic function, non-real zeros, unit disc. 
MSC 2010: 30D20, 30D35.
\end{abstract}
\section{Introduction} 

Around 1911 Wiman conjectured that if $f$ is a real entire function such that $f$
and $f''$ have only real zeros, then $f$ belongs to
the Laguerre-P\'olya class $\mathcal{LP}$
consisting of
locally uniform limits of real polynomials
with real zeros.
As the first main step towards the proof of this conjecture, it was shown in
 \cite{LeO} that for $f$  as in the assumptions the maximum modulus
$M(r, f)$ satisfies
\begin{equation}
\log^+ \log^+ M(r, f) = O( r \log r ) \quad \hbox{as $r \to + \infty$.} 
\label{LeOest}
\end{equation}
Wiman's conjecture was proved in full in \cite{BEL,SS}, while 
extensions to higher derivatives, as well as  to
real meromorphic functions (that is, meromorphic functions mapping
$\R$ into $\R \cup \{ \infty \}$), may be found in  
\cite{BEpolya,  HSW, Hin3,  lajda, Lajda18} and elsewhere. The aim of the present paper is to prove an estimate in the direction of (\ref{LeOest}) on the unit disc $\Delta =
D(0, 1)$, where
$D(a, r) =  \{ z \in \C : \, |z-a| < r \}$.

\begin{thm}
\label{thm1} 
Let the function $f$ be real meromorphic in the unit disc $\Delta $, that is,
 $f$ is meromorphic on $\Delta$ and maps the real interval $(-1, 1)$ into $\R \cup \{ \infty \}$. 
 Assume that $f^{(k)}$ has   finitely many non-real zeros in $\Delta$, for some $k \geq 2$, and that at least one of the following holds:\\
(i)   $f$ has finitely many poles and non-real zeros in $\Delta$;
\\
(ii) $1/f$ has finitely many poles and non-real zeros in $\Delta$.
\\
Then  $f$ satisfies, with
$h = f$ in case (i) and $h = 1/f$ in case (ii),
\begin{equation}
\label{f2est}
T(r, f)  \leq \log^+ M(r, h) + O(1)
\leq \exp \left( o(1-r)^{-4} \right) \quad \hbox{as $r \to 1-$. }
\end{equation} 
\end{thm}

Here $T(r, f)$ denotes the Nevanlinna characteristic function \cite{Hay2}, and
 it evidently suffices to prove the 
second inequality in (\ref{f2est}).  
Two key tools in the proof of the Wiman conjecture \cite{BEL,LeO,SS} were
Tsuji's analogue $\mathfrak{T}(r, f)$ for functions in a closed half-plane of
Nevanlinna's characteristic \cite{GO,Tsuji}, and  the Levin-Ostrovskii
factorisation $f'/f = P \psi $ of the 
 logarithmic derivative of a real entire function $f$ with only real zeros
 \cite{BEL,LeO}: here $P$ is a real entire function and either $\psi$ is constant or $\psi$ maps the open
 upper half-plane $H^+$ into itself.  In the proof of Theorem~\ref{thm1},  an analogue of the  Levin-Ostrovskii 
factorisation goes through relatively straightforwardly (see Lemma \ref{lemlevost}). 
However, as far as the author is aware, there is no direct counterpart of 
 the Tsuji characteristic for functions meromorphic in a half-disc. Some estimates are obtained
 nevertheless via a conformal mapping in Lemma \ref{lemtsuji1}, this requiring  slight modifications to the assumptions
 underlying the original Tsuji characteristic (see Propositions \ref{propfirsttsuji} and \ref{propsecondtsuji}).
Similar methods  lead to the following.

\begin{thm}
\label{thm2}
 Suppose that the function $f$ is meromorphic in $\Delta$, and assume that:\\ 
 (i) all but finitely many zeros and poles of $f$ in $\Delta$ are real;\\
 (ii)   
 for some $k \geq 0$, all but finitely many $1$-points of $f^{(k)}$ in $\Delta$ are real;\\ (iii) 
$f$ has a Nevanlinna deficient value, that is,
there exists $a \in \C \cup \{ \infty \}$ with
$$\delta (a, f) = \liminf_{r \to 1-}
\frac{m(r,a,f)}{T(r, f)}
> 0.$$
Then
\begin{equation}
T(r, f) \leq o(1-r)^{-3} \quad \hbox{as $r \to 1-$.}
\label{admissibleest}
\end{equation}
\end{thm} 

Note that the deficiency $\delta(a,f)$
is defined in \cite[pages 6, 42]{Hay2} only for functions which are admissible by virtue of growing fast enough to satisfy
\cite[formula (2.8), p.41]{Hay2},
but  non-admissible functions automatically satisfy (\ref{admissibleest}) anyway. Furthermore,
by a result of Edrei \cite[Theorem 1]{Edrei}
concerning roots on a system of rays, if a meromorphic function $f$ in the plane has only real zeros and poles, and
all $1$-points of $f^{(k)}$ are real for some $k \geq 0$, then $f$ has order $\rho(f) \leq 1$.


It seems very unlikely that either of Theorems \ref{thm1} and \ref{thm2}  is  sharp, and the exponents $-4$ and $ -3$,
arising in  conclusions (\ref{f2est}) and (\ref{admissibleest}) respectively, are almost certainly not best possible. However, a  simple example
shows that the second estimate of (\ref{f2est}) is
not too far out. Write
\begin{eqnarray}
w &=& \frac{1+z}{1-z} , \quad \frac{dw}{dz} = \frac2{(1-z)^2} ,
\quad f = e^{3e^w}, \quad
 L = \frac{f'}f = 3e^w w',
 \nonumber \\
L' &=&  3 e^w (w')^2 + 3 e^w w'' = 3 e^w  (w')^2 \left(  1 + \frac{w''}{(w')^2} \right) = 3 e^w (w')^2 \left(  2-z  \right)  ,
\nonumber \\
\frac{f''}{f} &=& L^2 + L' =
3 e^w (w')^2 \left( 3e^w + 2-z  \right) .
\label{notLP}
\end{eqnarray}
Then ${\rm Re} \, w > 0$ and $|e^w| > 1$ on  $\Delta$,
while
$
|2-z| < 3 < | 3  e^w | $ there.
This  makes $f$, $f'$ and  $f''$ all zero-free in the unit disc, and $h = f$ satisfies
$\log  M(r, h) = 3 \exp( (1+r)/(1-r))$.
Moreover,
(\ref{notLP})  shows that $(f'/f)' = L'$ is positive on the real interval $(-1, 1)$, whereas
if $P$ is a non-constant polynomial with only real zeros then a standard calculation shows that $(P'/P)' (x)$ is negative or infinite for each $x \in \R$. 
Thus $f$ is certainly not the locally uniform limit on $\Delta$ of a sequence of polynomials with only real zeros,
 and it does not seem obvious  what form
any direct counterpart of the Wiman conjecture should take in the unit disc setting. 

Furthermore, if $ g = 4 e^w$,
then $g$ belongs to the class of functions $G$ satisfying $ G(z) \neq 0, \infty$ and $G'(z) \neq 1 $ on $\Delta$,
for which
 growth estimates  may be found in \cite{SheaSons}.
Of course $g$ itself has bounded characteristic, since $|g(z)| > 4$ on $\Delta$, but
$\log M(r, g)$ grows like $(1-r)^{-1}$.
\hfill$\Box$
\vspace{.1in}

\section{The definition of the Tsuji characteristic} \label{tsujidef}

Let the function $f$ be meromorphic and not identically zero on the open upper half-plane $H^+$. 
For $r \geq 1$, set
\begin{equation}
\mathfrak{m}(r, f) =
 \frac1{2 \pi} \int_{ \arcsin(1/r)}^{\pi-\arcsin(1/r)} \frac{ \log^+ |f(r \sin \theta e^{i \theta} )|}{r \sin^2 \theta } \, d \theta ,
 \quad \mathfrak{N}(r, f) = \int_1^r \frac{  \mathfrak{n}(s, f) }{s^2} \, ds ,
 \label{ts5a}
\end{equation}
in which $\mathfrak{n}(r, f)$ is the number of poles of $f$, counting multiplicities,
in the set $\{ z \in \C :  \, |z| \geq 1 , \, |z-ir/2| \leq r/2 \}$. 
Tsuji's characteristic  function $\mathfrak{T}(r, f)$ is then given by \cite{GO,Tsuji} 
\begin{equation}
\mathfrak{T}(r, f) = \mathfrak{m}(r, f) + \mathfrak{N}(r, f) .
\label{ts7}
\end{equation}
The
key
properties associated with the Tsuji characteristic, including the fact  that it admits a lemma of the logarithmic derivative
\cite{GO,LeO,Tsuji},
 are normally stated and
proved when $f$ is meromorphic on the closed upper half-plane ${\rm Im} \, z \geq 0$:
this requirement is too strong for the present application, but the following two propositions will suffice. 

\begin{prop}
\label{propfirsttsuji} 
Suppose that the function $f$ is 
meromorphic and non-constant  on a domain containing the 
open upper half-plane $H^+$ and the points $1$ and $-1$. 
Then $\mathfrak{T}(r, f)$ may be defined by (\ref{ts5a}) and (\ref{ts7}), and 
differs from a continuous non-decreasing function by a bounded quantity. 
Moreover,  
 the first fundamental theorem holds for $f$, that is,  for any $a \in \C$,
  \begin{equation}
\label{ts12}
 \mathfrak{T}(r, f) = \mathfrak{T}(r, 1/(f-a)) + O(1) \quad \hbox{as $r \to + \infty $.} 
 \end{equation} 
 \end{prop} 
 It follows from (\ref{ts12}) and the fact that
  $1/z $ is bounded as $z \to \infty$ that the Tsuji characteristic of any rational function is bounded as $r \to + \infty$.
  
 \begin{prop}
 \label{propsecondtsuji}
  Suppose that the function $f$ is 
meromorphic and non-constant  on a domain containing $H^+$ and  the real interval $[-1, 1]$. 
Then  $f$ satisfies 
\begin{equation*}
 \mathfrak{m}(r, f'/f) \leq O( \log r + \log^+  \mathfrak{T}(r, f )  )
  \end{equation*}
 as  $r \to + \infty $ outside a set of finite linear measure. 
 \end{prop} 
 
 It is not claimed here that Propositions \ref{propfirsttsuji} and \ref{propsecondtsuji} are necessarily new, and relaxed conditions for the application of Tsuji's characteristic are referred to in various places,
 but not always with full detail (see e.g. \cite{LeO}). Proofs, which closely follow
  those in \cite{GO,Tsuji},
 will be given for completeness in Sections 
 \ref{tsujimoddef} and \ref{logderproof}. For now, it is convenient to note that 
 if 
$r \geq 1$ and $z$ lies in the set
\begin{equation}
\label{ts1}
J(r) = \{  r \sin \theta e^{i \theta } : \, 
\arcsin (1/r) \leq \theta \leq \pi - \arcsin (1/r)  \} ,
\end{equation}
which corresponds to the range of integration  in  $\mathfrak{m}(r, f)$, then $z$ satisfies
\begin{equation*}
|z| = \rho = r \sin \theta \geq 1, \quad 
| z - ir/2 |^2 = 
\rho^2 \cos^2 \theta + \rho^2 \sin^2 \theta - r \rho \sin \theta + \frac{r^2}4 =  \frac{r^2}4 .
\end{equation*}
Thus $J(r)$ is the arc of the circle $|z-ir/2| = r/2$ going counter-clockwise from 
$e^{ i \arcsin (1/r) }$ to
$e^{ i ( \pi - \arcsin (1/r) )}$
via $ir$, and
(see \cite{LeO} and \cite[p.322]{GO})
\begin{equation}
\{ \zeta \in \C: \, |\zeta| \geq 1, 
\, {\rm Im} \, \zeta  > 0 \} \subseteq \bigcup_{r \geq 1} J(r).
\label{J(r)inclusion}
\end{equation}


\section{The Levin-Ostrovskii factorisation in the unit disc}

\begin{lem} 
\label{lemlevost}
Let the function $g$ be real meromorphic in  $\Delta$ such that $g$ has 
 finitely many poles and non-real zeros. Then $g'/g$ has a Levin-Ostrovskii factorisation 
$g'/g = P \psi$, in which $P$ and $\psi$ are real meromorphic functions on $\Delta$  such that: (i)
$P$ has finitely many poles;  (ii)
either $\psi \equiv 1$ or $\psi $ maps the upper half-disc 
$\Delta^+ = \{ z \in \C : \, |z| < 1, \, {\rm Im} \, z > 0 \} $
into the upper half-plane $H^+ = \{ z \in \C : \, {\rm Im} \, z > 0 \}$.
\end{lem}
\textit{Proof.} If $g$ has finitely many zeros, set $\psi \equiv 1$. Assume henceforth that $g$ has infinitely many 
zeros: then all but finitely many of the real zeros of $g$  may be labelled $a_k$ in such a way that Rolle's theorem gives real zeros $b_k$ of $g'$ satisfying
$
a_k < b_k < a_{k+1} ,
$
for $-\infty \leq k_1 < k < k_2 \leq + \infty$, say.
If $|z| \leq r < 1$ and 
$|k|$ is large enough, then
$$
\left| 1 - \frac{b_k-z}{a_k-z} \right| = \left| \frac{a_k-b_k}{z-a_k} \right|
\leq \frac{2 (b_k - a_k)}{1-r}  .
$$
Since
$\sum_{k_1 < k < k_2} (b_k-a_k) < 2$, 
the product
$$
\psi (z) = \prod
_{k_1 < k < k_2}
\left( \frac{b_k-z}{a_k-z} \right)
$$
converges uniformly on compact subsets of $\Delta$ and satisfies, for $z \in \Delta^+$, 
$$
0 < \arg \psi (z) = \sum
_{k_1 < k < k_2} ( \arg (b_k-z) - \arg (a_k-z) ) < \pi .
$$
\hfill$\Box$
\vspace{.1in}

\section{A growth estimate  on the unit disc}

\begin{lem}
\label{lemwmap}
The function 
\begin{equation} 
\label{wdef}
w = w(z) = \frac{4z}{z^2+1} 
\end{equation}
maps: (a) $\Delta$ univalently onto $\C \setminus \{ X \in \R : \, |X| \geq 2 \}$;   (b)
$\Delta^+$  onto $H^+$; (c)
 the  real interval $(-1, 1)$ onto $(-2, 2)$; (d)
$[\sqrt{3}-2, 2 - \sqrt{3}]$ onto $[-1, 1]$. Moreover, $|w| \geq 1$ for
$ 2- \sqrt{3} \leq |z| < 1$.
\end{lem}
\textit{Proof.} 
The well known mapping $p(z) = z+1/z = 4/w$
satisfies $p(u) = p(v)$ if and only if $u=v$ or $u = 1/v$. Furthermore, $p$ maps the unit circle $|z| = 1$ onto the interval
$[-2, 2]$ and $\Delta^+$ onto the open lower half-plane,
while 
$w(x)$ is real and increasing for   $x \in (-1, 1)$.
The last assertion holds since
$|w(z)| \geq w(|z|)$.
\hfill$\Box$
\vspace{.1in}

\begin{lem}
\label{lempsidisc}
Suppose that the function  $\psi$ maps
$\Delta^+$ analytically into $H^+$.
Then  $\psi$ satisfies,
for $z \in \Delta^+$ with $ 2- \sqrt{3} \leq |z| < 1$, 
$$
\frac{  \left(1-|z|^2 \right) { \rm Im} \, z}{20 |z|^2 } \leq \left| \frac{\psi(z)}{\psi((\sqrt{5}-2)i)} \right| \leq 
\frac{20 |z|^2 }{  \left(1-|z|^2 \right) { \rm Im} \, z}.
$$
\end{lem}
\textit{Proof.} Define $w$  by (\ref{wdef}) and set
$ \Psi(w) = \psi(z)$.
Then $\Psi : H^+ \to H^+ $ is analytic and so \cite{BEL,Lev}
$$
\frac{{\rm Im} \, w }{5|w|^2}  \leq \left| \frac{\Psi(w)}{\Psi(i) } \right| \leq \frac{5|w|^2}{{\rm Im} \, w }  \quad \hbox{
for $w \in H^+$ with $|w| \geq 1$.
}
$$
The asserted inequalities
then follow from Lemma \ref{lemwmap}, since
$w((\sqrt{5}-2)i) = i$ and
$$
\frac{2i \,  {\rm Im} \, w }{|w|^2} = \frac{4}{|w|^2} \left(\frac{z}{1+z^2} -  \frac{\overline{z}}{1+ \overline{z}^2 } \right) 
= \frac{4}{16 |z|^2} \left(z \left( 1+ \overline{z}^2 \right) -  \overline{z} \left(1+z^2 \right)  \right) =
\frac{2i  \left(1-|z|^2 \right) {\rm Im} \, z }{4 |z|^2 } .$$
\hfill$\Box$
\vspace{.1in}



\begin{lem}
\label{lemrational} 
Suppose that $R_1(w) = R_2(z)$ for  $z \in \Delta$ and  $w \in \C \setminus \{ X \in \R : \, |X| \geq 2 \}$ as in (\ref{wdef}),
where $R_2$ is a rational function. Then
$\mathfrak{T}(r, R_1)$ is bounded as $r \to + \infty$.
\end{lem}
\textit{Proof.} 
The assumptions imply that  $R_2(z) $ is a quotient of bounded functions on $\Delta$,  
and therefore so is $R_1(w)$ on $\C \setminus \{ X \in \R : \, |X| \geq 2 \}$, by Lemma \ref{lemwmap}.
Hence the conclusion follows from  Proposition  \ref{propfirsttsuji} and standard properties of the Tsuji characteristic.
\hfill$\Box$
\vspace{.1in}
   

\begin{lem}
\label{lemtsuji1}
Let  $f$ be a meromorphic function on the unit disc $\Delta$, and write $g(w) = f(z)$,
with $w$ as in  (\ref{wdef}). 
Then $g$ is meromorphic on $\C \setminus  \{ X \in \R : \, |X| \geq 2 \}$.
Moreover, if 
\begin{equation}
\label{tsuji0}
I_1 = 
\int_1^{+ \infty} \frac{\mathfrak{m}(t,g)}{t^2} \, dt < + \infty 
\end{equation}
 then
\begin{equation}
\label{tsuji2}
\int_{2-\sqrt{3}}^1 \int_0^\pi (1-r^2 )^2 \log^+ |f(re^{i \theta} ) | \, d \theta \, dr < + \infty .
\end{equation}
In particular, (\ref{tsuji2}) holds if 
\begin{equation}
\label{tsuji1}
\mathfrak{T}(t, g) = O( \log t )  \quad \hbox{ as $t \to + \infty $.}
\end{equation}
\end{lem} 
\textit{Proof.} 
The first assertion follows from Lemma \ref{lemwmap}.
Assume now that (\ref{tsuji0}) holds, which will certainly be the case if  (\ref{tsuji1}) is true, and write 
$s = t \sin \theta $.
Then (\ref{J(r)inclusion})  yields, following \cite{LeO} and \cite[Ch. 6, Lemma 5.2]{GO},
\begin{eqnarray*}
+\infty > I_1 
&=&
 \frac1{2 \pi} \int_1^{+\infty} \int_{\arcsin(1/t)}^{\pi - \arcsin(1/t) } 
\frac{\log^+ | g(t \sin \theta e^{i \theta })| }{t^3 \sin^3 \theta } \, \sin \theta \,  d \theta \, dt \\
&\geq&
 \frac1{2 \pi} \int_{1 \leq s < + \infty, \, 0 < \theta < \pi } 
\frac{\log^+ | g(s e^{i \theta })| }{s^3 } \,  d \theta \, ds = I_2 .
\end{eqnarray*}
If $z =r e^{i \phi} $ with $2-\sqrt{3} \leq r < 1$ and $0 < \phi < \pi $ then $w \in H^+$ and 
$|w| \geq 1$, by Lemma \ref{lemwmap}.
Hence 
writing $z = x+iy$ and $w = u+iv = se^{i \theta} $ with $x, y, u, v$ all real  leads to 
\begin{eqnarray*}
+\infty &>& I_2 = 
\int_{|w| \geq 1, v > 0} \frac{\log^+ |g(w)|}{|w|^4 } \, du \, dv 
\geq
\int_{1 > |z| \geq 2-\sqrt{3} , \, y > 0} \frac{\log^+ |f(z)|}{|w|^4 } \left| \frac{dw}{dz} \right|^2 \, dx \, dy .
\end{eqnarray*}
This delivers (\ref{tsuji2}), in view of the fact that, for $|z| < 1$,
$$
\frac1{|w|^4}\left| \frac{dw}{dz} \right|^2 = \left| \frac{(1-z^2)^2}{16 z^4 } \right| \geq \frac{(1-|z|^2 )^2}{16} .
$$
\hfill$\Box$
\vspace{.1in}

\section{Proof of Theorem \ref{thm2}} 

Let $f$ be as in the hypotheses  of Theorem \ref{thm2}.
It may be assumed that $f$ is not  rational, and that $f$ is admissible in the sense of \cite[p.42]{Hay2}, as otherwise the conclusion of the theorem is obvious.

It will be shown that  (\ref{tsuji1}) holds for
the function $g$ defined by $g(w) = f(z)$ and (\ref{wdef}), which
is meromorphic on
$\C \setminus  \{ X \in \R : \, |X| \geq 2 \}$ by Lemma \ref{lemtsuji1}. Assume that this is not the case, so that
$\mathfrak{T}(r, g)$ must be unbounded and tend to infinity, and 
 observe that 
assumption (i) and  Lemma \ref{lemwmap} imply that $g$ has finitely many zeros and poles in $H^+$. Moreover, 
there exist rational functions $b_0(z) , \ldots, b_k(z)$ such that 
$$
f^{(k)}(z) = \sum_{j=0}^k b_j(z) g^{(j)}(w) =  \sum_{j=0}^k c_j(w) g^{(j)}(w) = \phi(w), \quad
c_j(w) = b_j(z),
$$
in which each $c_j$ has bounded Tsuji characteristic by Lemma \ref{lemrational}. Thus $\phi$ has finitely many 
$1$-points in $H^+$, and is non-constant, since otherwise $f$ is a polynomial, contrary to assumption. 
Hence Milloux' theorem \cite[p.57]{Hay2} (when $k \geq 1$) or the second fundamental theorem 
(when $k=0$) may be applied, with the Tsuji characteristic in place of that of Nevanlinna,
and so $g$ satisfies (\ref{tsuji1}) as asserted.
Furthermore, if $F(z) = \bar f (\bar z) $ then $F^{(k)}-1$ has finitely many non-real zeros in $\Delta$,
and $G(w) = F(z)$ also has
$\mathfrak{T}(r, G) = O( \log r )$.

 Suppose first
 that $a = \infty$ in assumption (ii), and set $h=f$ and
$$
I_3 =
\int_{2-\sqrt{3}}^1 \int_{- \pi}^{\pi} (1-r^2 )^2 \log^+ |h(re^{i \theta} ) | \, d \theta \, dr .$$
 Then
Lemma \ref{lemtsuji1}
applied to $f$ and $F$ delivers
$I_3 < + \infty$ and, as $r \to 1-$,
since $\delta (\infty, h) > 0$,
$$
T(r, h) \int_r^1 (1-t )^2  \, dt \leq \int_r^1  T(t, h) (1-t )^2  \, dt \to 0, \quad
  T(r, h)   = o(1-r)^{-3},$$
as required.
Now suppose that
$a \in \C$, and set
$h = 1/(f-a)$. Then
$$
T(r, f) = T(r, h) + O(1), \quad
\mathfrak{T}
(t, 1/(g-a)) +
\mathfrak{T}
(t,1/(G - \bar a) ) = O( \log t ),
$$
as $r \to 1-$ and $t \to + \infty$,
by the first fundamental theorem,
and this time
\begin{eqnarray*}
 I_3 &=&
 \int_{2-\sqrt{3}}^1 \int_{0}^{\pi} (1-r^2 )^2 \log^+ |1/(f(re^{i \theta} ) -a ) | \, d \theta \, dr
 + \\
 & & + \int_{2-\sqrt{3}}^1
 \int_{0}^{\pi} (1-r^2 )^2 \log^+ |1/(F(re^{i \theta} ) - \bar a ) | \, d \theta \, dr < + \infty ,
\end{eqnarray*}
by Lemma \ref{lemtsuji1},
which gives
$ T(r, h)   = o(1-r)^{-3}$
as before.
 \hfill$\Box$
\vspace{.1in}

\section{ Wiman-Valiron theory in the unit disc}\label{wimanvaliron}

Let the function $h$ be meromorphic in $\Delta= D(0, 1)$,
and assume that $h$ has  
finitely many poles in $\Delta$, and that $|h(z)| $ is unbounded as $|z| \to 1-$. Choose $s_0 \in (0, 1)$ such that all poles of $h$ in
$\Delta$ lie in the disc $D(0, s_0)$,  and take $R > 0$ such that $|h(z)| < R$ on
$|z| = s_0$.
The set 
\begin{equation}
U = \{ z \in \C : \, s_0 < |z| < 1, \, |h(z)| > R \} 
\label{setUdef}
\end{equation}
is then non-empty and open, and it follows from the maximum principle that if $C$ is a component of $U$ 
then  $\partial C$ meets the circle  $|z| = 1$. 
Next, the function 
\begin{equation}
 \label{1}
v(z) = \log \left| \frac{h(z)}{R} \right| \quad (z \in U), \quad v(z) = 0 \quad (z \in \Delta\setminus U) ,
\end{equation}
is  continuous, subharmonic and unbounded on $\Delta$. 
For $0 < r < 1$ let 
\begin{equation}
 \label{2}
B(r) = B(r, v) = \max \{ v(z) : |z| = r \} , \quad a(r) = r B'(r) = \frac{dB(r)}{d \log r } .
\end{equation}
Here $B(r)$ is a non-decreasing convex function of $\log r$ for $0 < r < 1$. Moreover, 
$a(r)$ (which is  taken to be the right derivative  at the countably many 
points at which $B$ is not differentiable) is non-decreasing and tends to $+ \infty$ as $r \to 1-$.
The following result of Wiman-Valiron type \cite{Hay5} is a modified version of \cite[Theorem 1.1]{LaRo2013}.

\begin{prop}
 \label{propA}
Choose $r_0 \in (0, 1)$ with $B(r) \geq 2$ and $a(r) \geq 2 $ for $r_0 \leq r < 1$, as well as
$\beta \in (0 , 1/2]$ and $  \delta > 0$, and  set
\begin{equation}
 \label{3a}
\varepsilon (r) = \min \left\{ \frac{1-r}{2 a(r)^\beta (\log a(r))^{1 +  \delta } }, \quad  
\frac{ 1}{  a(r)^{1-\beta}(\log a(r))^{1 +  \delta } }\right\} 
\end{equation}
for $r_0 \leq r < 1$. Then there exists a set $E \subseteq [r_0, 1)$ satisfying
\begin{equation}
 \label{4}
\int_E \, \frac{dt}{1-t} < \infty ,
\end{equation}
such that, as $r \to 1-$ with $r \not \in  E$, if $z_r$ is chosen with
$|z_r| = r$ and $v(z_r) = B(r, v)$ then 
\begin{equation}
 \label{con1}
h(z) \sim h(z_r) \left( \frac{z}{z_r} \right)^{a(r)} 
 \quad \hbox{and} \quad 
\frac{h'(z)}{h(z)} \sim \frac{ a(r)}z  \quad \hbox{ for} \quad  | z- z_r | < \frac{ \varepsilon (r)}{2048}.
\end{equation}
\end{prop}

The proof of Proposition \ref{propA} is almost identical to that of 
 \cite{LaRo2013}, which was in turn adapted from that of \cite{BRS}, but
 one observation should be made: the set $U$ in (\ref{setUdef}) need not be connected, whereas
 \cite{LaRo2013} considered  meromorphic functions 
on $\Delta$ with a direct tract $U$ in the sense of \cite{BRS}, in which case $U$ is connected by definition. However,
this does not affect the proof, some key ingredients  of which may be outlined  as follows.

The set $E$ is determined by the following lemma, a combination of  \cite[Lemmas 3.1 and 3.2]{LaRo2013}:
the proof, which is adapted in \cite{LaRo2013}  from \cite{bergsizediscs},
 depends only on (\ref{2}) and the fact  that $a(r)$ is non-decreasing and
unbounded on $(r_0, 1)$.

\begin{lem}
 \label{lemexset}
There exists a set $E \subseteq [r_0, 1)$ satisfying (\ref{4}) 
such that, for $r $ in $ [r_0, 1) \setminus E$, 
\begin{equation*}
 \label{7}
a(r + \varepsilon(r)) < a(r) + a(r)^{1-\beta} , \quad a(r - \varepsilon(r)) > a(r) - a(r)^{1-\beta} ,
\quad (1-r) a(r) < B(r)^{1 + \beta } .
\end{equation*}
Moreover, $B(r)$ satisfies, for $r \in [r_0, 1) \setminus E$, 
\begin{equation}
 \label{14}
B(s) \leq B(r) + a(r) \log \frac{s}r + \phi(r)  \quad \hbox{for} \quad 
r - \varepsilon (r) \leq s \leq r + \varepsilon (r) ,
\end{equation}
in which $0 \leq \phi(r)
\leq O( a(r)^{1-\beta} \varepsilon(r) )
\to 0$  as $r \to 1-$ in $ [r_0, 1) \setminus E$.
\end{lem}
\hfill$\Box$
\vspace{.1in}

\begin{lem}
 \label{disclem}
Let $r \in [r_0, 1) \setminus E$, set $\sigma = \sigma(r) = \varepsilon (r)/2048$ and choose $z_r $
with $|z_r| = r$ and $v(z_r) = B(r, v)$. If $r $ is close enough to $1$ then the disc $D(z_r, 4 \sigma) $  lies in $ U$. 
\end{lem}
As in  \cite[Lemma 3.3]{LaRo2013}, the proof begins by noting that, 
on $D(z_r, 2048 \sigma)$, by 
(\ref{14}),
\begin{equation}
 \label{16}
u(z) = v(z) - B(r) - a(r) \log \frac{|z|}r \leq \phi(r) = o(1) .
\end{equation}
Take $R'$ with $R' - R$ small and positive such that $h$ has no critical values of modulus $R'$,
and 
assume that $r \in [r_0, 1) \setminus E$ is close to $1$ and that
$D(z_r, 4 \sigma) \not \subseteq  U$.
Let $U'$ be that component of $\{ z \in \C : \, s_0 < |z| < 1, \, |h(z)| > R' \} $ which contains $z_r$:
then there is a component $K$ of $\Delta \setminus U'$ which meets 
$D(z_r, 4 \sigma) $. A contradiction is then obtained exactly as in \cite{LaRo2013} and, because $U'$ is
automatically connected, it makes no difference that $U$ might not be. 
\hfill$\Box$
\vspace{.1in}

To establish Proposition \ref{propA}, take  $r \not \in E$, close to $1$, and
apply  the Borel-Carath\'eodory inequality 
to  (\ref{16}), as in  
 \cite[Lemmas 3.4 and 3.5]{LaRo2013}, to obtain the estimate for $h(z)$ in (\ref{con1}) via
\begin{equation*}
 \label{19}
\log \frac{h(z)}{h(z_r)} = a(r) \log \frac{z}{z_r} + g(z), \quad | g(z)| \leq 2 \phi (r) = o(1)
\quad \hbox{ on $z \in D(z_r, 2 \sigma)$.}
\end{equation*}
Cauchy's estimate for derivatives  then gives the estimate for $h'(z)/h(z)$
on $ D(z_r,   \sigma)$.
\hfill$\Box$
\vspace{.1in}



\section{The Tumura-Clunie method on a half-disc} 
\label{clunie}

Let $G $ be a non-constant
meromorphic function on the unit disc $\Delta$. Using the same device as in Lemma~\ref{lemtsuji1},  define $G^*$ by
$G^*(w) = G(z)$, with $w$ as in (\ref{wdef}). Then $G^*$ is meromorphic on $\C \setminus \{ X \in \R : \, |X| \geq 2 \}$,
by Lemma \ref{lemwmap},
and so has a Tsuji characteristic: for $r \geq 1$ write
\begin{equation*}
\mathfrak{m}^*(r, G) = \mathfrak{m}(r, G^*), \quad
\mathfrak{N}^*(r, G) = \mathfrak{N}(r, G^*), \quad
\mathfrak{T}^*(r, G) = \mathfrak{T}(r, G^*) .
\end{equation*} 
Denote by
$\mathfrak{S}^*(r, G)= \mathfrak{S}(r, G^*)
$  any term which satisfies
$$\mathfrak{S}(r, G^*) \leq
O(  \log r + \log^+ \mathfrak{T}(r, G^*) )
 $$
as $r \to +\infty$ outside a set of finite measure.
Thus Propositions
\ref{propfirsttsuji} and
\ref{propsecondtsuji} imply that
the first fundamental theorem
and the lemma of the logarithmic derivative
hold for   $\mathfrak{T}^*$,  since Lemma \ref{lemrational} and the formula
\begin{equation}
G'(z) = (G^*)'(w) \frac{dw}{dz} =  (G^*)'(w) \frac{4(1-z^2)}{(1+z^2)^2}
\label{Fdiff}
\end{equation}
together give
$
\mathfrak{m}^*(r, G'/G) = \mathfrak{S}^*(r, G)
$. Let
$\Lambda_G$ be the set of meromorphic functions
$H$ on $\Delta$
satisfying $\mathfrak{T}^*(r, H) = \mathfrak{S}^*(r, G)$.
As is standard, a differential polynomial in $G$ over $\Lambda_G$
means a polynomial in $G$ and its derivatives, with coefficients in $\Lambda_G$, and a result
of Tumura-Clunie type \cite[p.67]{Hay2} goes through easily as follows.

\begin{lem}
\label{lemclunie} 
Assume that 
$$
G(z)^n P(z) = Q(z),
$$
where $P$ and $Q$ are differential polynomials in $G$ over $\Lambda_G$. If $Q$ has
total degree at most $n$ in $G$ and its derivatives, then $\mathfrak{m}^*(r, P) = \mathfrak{S}^*(r, G)$ as $r \to + \infty$.

Furthermore, if  $H_{k-1}(G) $ is a  differential polynomial in $G$ over $\Lambda_G$, of total degree at most $k-1$, and if
$$
\Phi = G^k + H_{k-1}(G),
\quad
\mathfrak{N}^*(r, G) + \mathfrak{N}^* (r, 1/\Phi ) = \mathfrak{S}^*(r, G),$$
then 
either $\Phi = (G + \alpha)^k $ for some $\alpha \in \Lambda_G$,
or $G \in \Lambda_G$ and $\mathfrak{T}^*(r, G) = O(\log r)$  as $r \to + \infty$.
\end{lem} 
\textit{Proof.} 
To prove the first part it  suffices to write the equation as 
$
G^*(w)^n P^*(w) = Q^*(w)
$
and observe that, by Lemma \ref{lemrational} and repeated differentiation of (\ref{Fdiff}),
$P^*(w)$ and $Q^*(w)$ can be written as differential polynomials in $G^*(w)$,
their coefficients $a(w)$
satisfying
$\mathfrak{T}(r, a) =
 \mathfrak{S}(r, G^*) $.
The standard
proof of Clunie's lemma from \cite[p.68]{Hay2} may then be applied.

Next, it is clear  that if $G \in \Lambda_G$ then $\mathfrak{T}^*(r, G) = O(\log r)$  as $r \to + \infty$, so assume that
$G \not  \in \Lambda_G$. Then the conclusion $\Phi = (G + \alpha)^k $ for some $\alpha \in \Lambda_G$
follows from applying the standard proof more or less verbatim as in 
\cite[pp.69-73]{Hay2}. 
 \hfill$\Box$
\vspace{.1in}

The points  $w \in H^+$ which contribute to
$\mathfrak{m}(r, G^*) = \mathfrak{m}^*(r, G)$
and the associated $z \in \Delta^+$
are discussed at the end of Section \ref{tsujimoddef}.


\section{Proof of Theorem \ref{thm1}} 

Let $f$ be as in the statement. It may be assumed that $L = f'/f$ is not a rational function, since otherwise the conclusion of the theorem follows immediately.
Let $h$  be $f$ in case  (i), with  $h = 1/f$ in case (ii),
so that $h$ has finitely many poles and non-real zeros in $\Delta$.
Applying Lemma~\ref{lemlevost} to $h$ gives a factorisation
\begin{equation}
\label{Ldef} 
L = \frac{f'}{f} = P \psi ,
\end{equation} 
in which $P$ and $\psi$ are real meromorphic in $\Delta$,  such that $P$ has finitely many poles and 
either $\psi \equiv 1$ or $\psi ( \Delta^+) \subseteq H^+$. 
With $w$ as in (\ref{wdef}), set
\begin{equation} 
\label{MQdef}
L(z) = M(w), \quad P(z) = Q(w) .
\end{equation}
Then
$M$ and $Q$ are meromorphic on $\C \setminus \{ X \in \R : \, |X| \geq 2 \}$, by Lemma \ref{lemwmap}.

\begin{lem} 
\label{lemMestb}
The functions $L$ and $M$ satisfy $\mathfrak{T}^*(r, L) = \mathfrak{T}(r, M) = O( \log r )$  as $r \to + \infty $.
\end{lem} 
\textit{Proof.}  
Assume that this is not the case, so that $L \not \in \Lambda_L$. Lemma \ref{lemclunie} implies that, since $f^{(k)}/f$
is a differential polynomial in $L$ \cite[Lemma 3.5]{Hay2} with finitely many zeros and poles in $\Delta^+$,
\begin{equation}
H = \frac{f^{(k)}}f = \left( L + \alpha \right)^k , 
\label{TT1}
\end{equation}
where $\alpha \in \Lambda_L$ and so $L+\alpha \not \equiv 0$. For $j=1, \ldots, k$ set $f_j(z) = z^{j-1} $. Then 
 (\ref{TT1}) and the equation $f^{(k)} = W(f_1, \ldots, f_k, f)$ yield
\begin{eqnarray*}
\frac{\left( L + \alpha \right)^k}{f^k}  &=&   
 \frac{W(f_1, \ldots, f_k, f)} {f^{k+1}} =
 W(f_1/f, \ldots, f_k/f, 1)  
 = (-1)^k 
 W((f_1/f)', \ldots, (f_k/f)') ,
\end{eqnarray*}
and so 
\begin{equation}
 \label{gjdef}
W(g_1, \ldots, g_k) = (-1)^k, \quad
g_j = \left( \frac{f_j}f \right)' \, \frac{f}{L+ \alpha} = \frac{f_j' - L f_j}{L + \alpha } .
\end{equation}
By (\ref{TT1}), the $g_j$, which  are the same  auxiliary functions as  in Frank's method \cite{FHP},
satisfy
\begin{equation}
\mathfrak{T}^*(r, g_j) \leq O( \mathfrak{T}^*(r, L) )  + \mathfrak{S}^*(r, L) \quad \hbox{ as $r \to \infty$,}
\label{gjest}
\end{equation}
and each has finitely many poles in $\Delta^+$. Hence the $g_j$ are solutions on $\Delta$
of an equation
\begin{equation}
w^{(k)} + \sum_{j=0}^{k-2} a_j w^{(j)} = 0, \quad a_j \in \Lambda_L.
\label{TT2}
\end{equation}
If $g_p + f_p $ vanishes identically, then
$ -f_p \alpha = f_p' $, which cannot be the case for all $p$,
since $k \geq 2$.
Hence there exists at least
one $q$ with $g_q + f_q
\not \equiv 0$, in which case (\ref{gjdef}) yields
\begin{equation}
L = \frac{f_q' - \alpha g_q}{g _q+f_q}
= \frac{f_q' + \alpha f_q}{g _q+f_q} - \alpha .
\label{TT3}
\end{equation}
If $a_j \equiv 0$ for all $j$ then the $g_j$ are polynomials, and so
(\ref{TT3}) gives  $L \in \Lambda_L$, a contradiction.

Thus  $a_j \not \equiv 0$
for at least one $j$, and so at least one $f_q$
does not solve (\ref{TT2}), and in particular $g = g_q+f_q \not \equiv 0$.
Since $g_q$ solves (\ref{TT2}) but $f_q$ does not,
there exists $S \in \Lambda_L $ with
$$
g^{(k)} +  \sum_{j=0}^{k-2} a_j g^{(j)} = f_q^{(k)} +  \sum_{j=0}^{k-2} a_j f_q^{(j)} = S\not \equiv 0, \quad 
\frac1g = \frac1S \left( \frac{g^{(k)}}g  +  \sum_{j=0}^{k-2} a_j \frac{g^{(j)}}g  \right).
$$
On combination with
(\ref{gjest}),
(\ref{TT3}) and the fact that
 $L$ has finitely many poles in $\Delta^+$, this leads to
$$
\mathfrak{T}^*(r, L) \leq
\mathfrak{m}^*(r, L) + \mathfrak{S}^*(r, L) \leq \mathfrak{m}^*(r, 1/g) + \mathfrak{S}^*(r, L) \leq \mathfrak{S}^*(r, L),
$$
a contradiction.
\hfill$\Box$
\vspace{.1in}


\begin{lem} 
\label{lemQest}
The  function $Q$ in (\ref{MQdef})
satisfies $\mathfrak{T}(r, Q) = O( \log r )$ as $r \to + \infty $.
\end{lem} 
\textit{Proof.} Write $\Psi(w) = \psi(z)$, with $\psi$ and $w$ as in (\ref{Ldef}) and (\ref{MQdef}), so that 
$M = Q \Psi$ on $H^+$. The assertion then follows
from  Lemma \ref{lemMestb},  standard properties of the Tsuji characteristic and the existence of a M\"obius transformation
$T_1$ such that $T_1 \circ \Psi$ is bounded on $H^+$. 
\hfill$\Box$
\vspace{.1in}

The
 real meromorphic function $P$ on $\Delta$ has finitely many poles. Thus
 Lemmas \ref{lemtsuji1} and \ref{lemQest} imply that,
 with  $C$ denoting positive constants which depend at most on $L$ and $P$,
\begin{eqnarray*} 
\int_{2-\sqrt{3}}^1 (1-r)^2 T(r, P) \, dr &<& \int_{2-\sqrt{3}}^1 (1-r)^2 m(r, P) \, dr + C\\
&<& \frac1\pi 
\int_{2-\sqrt{3}}^1 \int_0^\pi (1-r)^2 \log^+ |P(re^{i \theta} ) | \, d \theta dr + C < C.
\end{eqnarray*} 
Since $T(r, P)$ is non-decreasing this yields, as $s \to 1-$, 
$$
T(s, P) \int_{s}^1 (1-r)^2  \, dr \leq
\int_{s}^1 (1-r)^2 T(r, P) \, dr  \to 0,
\quad T(s, P) = o( 1-s)^{-3} ,
$$
and, again since $P$ has finitely many poles, it follows that 
\begin{equation}
\log^+ M(r, P) = o( 1-r)^{-4} \quad \hbox{as $r \to 1-$.} 
\label{Pestimate}
\end{equation}

If $h(z)$ is bounded as $|z| \to 1-$ then (\ref{f2est})
evidently holds, so assume that $h(z)$ is unbounded
there and apply
Proposition \ref{propA},
with its notation.
As $r \to 1-$ with $r \not \in E$,  this  yields
\begin{equation}
\label{arcest}
|L(z)| = 
\left| \frac{f'(z)}{f(z)} \right| = 
\left| \frac{h'(z)}{h(z)} \right| \geq C a(r) 
\end{equation}
on an arc of the circle $|z| = r$ with length $C \varepsilon (r)$. Since
$0 < \beta \leq 1/2$ and $\psi$ is real meromorphic in $\Delta$,
Lemma \ref{lempsidisc},  (\ref{3a}), 
(\ref{Ldef}),  (\ref{Pestimate}) and (\ref{arcest}) deliver,
as $r \to 1-$ outside $E$,
\begin{eqnarray*}
a(r)^{1-\beta/2}  &\leq&
\frac{C}{\varepsilon (r)}
\int_0^{2\pi} | L(re^{i \theta } )|^{1-\beta/2} \, d \theta \leq
\frac{\exp \left( o(1-r)^{-4} \right)}
{\varepsilon (r)}
\int_0^{2\pi} | \psi(re^{i \theta } )|^{1-\beta/2} \, d \theta
\nonumber \\
  &\leq&
 \exp \left( o(1-r)^{-4} \right) 
\max \left\{ \frac{2a(r)^\beta (\log a(r))^{1+\delta } }{1-r} , \quad 
a(r)^{1-\beta} (\log a(r))^{1+\delta } \right\} ,
\end{eqnarray*}
and hence
$a(r) \leq \exp \left( o(1-r)^{-4} \right)$.
Moreover, the last estimate
 holds for $r$ close to $1$,  with no   exceptional set,
since $a(r)$ is non-decreasing and   $[r, (1+r)/2] \not \subseteq E$, by
(\ref{4}).
Integration  then yields
$
B(r) = B(r, v)
\leq \exp \left( o(1-r)^{-4} \right)
$
as $r \to 1-$,
and hence (\ref{f2est}).

\hfill$\Box$
\vspace{.1in}

\section{Proof of Proposition \ref{propfirsttsuji}} \label{tsujimoddef}

As already remarked, it is not claimed that Propositions  \ref{propfirsttsuji} and
 \ref{propsecondtsuji} are new, and the following proofs are included mainly for completeness. 
 Suppose  that  the function $f $ is 
meromorphic and non-constant
on a domain containing $H^+ \cup \{ -1, 1 \}$. Then so is the function $F$ defined by 
\begin{equation}
\label{ts3}
F(\zeta) = f(z) , \quad 
\zeta = - \,  \frac1z = t + i \sigma  ,
 \end{equation}
in which  $t, \sigma$ are real. 
This section follows closely the method from \cite{Tsuji}, 
but a slight simplification arises from choosing  $x \in \R$, close to $1$,
such that $F$ has no zeros or poles on the stepwise curve $L_x$
from $x$ to $-x$ via $x+ix$ and $-x+ix$. For $0 < \sigma \leq x$ let $K(\sigma) $ be the rectangle 
\begin{equation}
\label{ts9}
K(\sigma) = \{ u+iv : \, -x \leq u \leq x, \, \sigma \leq v \leq x \} .
\end{equation}
If $g$ is a  meromorphic function on the compact plane set $X$, denote by $n( X, g)$  the number of poles of $g$ in $X$, counting multiplicities, and  for $0 < \sigma \leq x$ set
\begin{eqnarray}
\label{ts10}
T_2(\sigma, F) &=& m_2(\sigma, F) + N_2(\sigma, F) =
\frac1{2 \pi } \int_{-x}^x \log^+ | F(t + i \sigma )| \, dt + N_2(\sigma, F), \nonumber \\
n_2(\sigma, F) &=& n(K(\sigma), F), \quad 
N_2(\sigma, F) = \int_\sigma^x n_2(\lambda, F) \, d \lambda.
\end{eqnarray}
\begin{lem}
\label{lembounded}
As $ \sigma \to 0+$ and $r = 1/\sigma \to + \infty $, 
\begin{equation}
m_2(\sigma, F) = \mathfrak{m}(r, f) + O(1) , \quad N_2(\sigma, F) = \mathfrak{N}(r, f) + O(1),
\quad
T_2(\sigma, F) = \mathfrak{T}(r, f) + O(1).
\label{101formula}
\end{equation}
Moreover, these relations also hold with $F$ and $f$ replaced by their reciprocals. 
\end{lem} 
\textit{Proof.} 
Under the change of variables (\ref{ts3}) 
 the circular arc $J(r)$  in (\ref{ts1}) transforms to the horizontal line segment $I(r)$  given by 
\begin{equation}
\label{ts4}
\zeta = - \, \frac{ \cot\theta}r + \frac{i}r , 
\quad \sigma = \frac1r, 
\quad t =  - \, \frac{ \cot \theta}r, \quad \frac{dt}{d\theta} = \frac1{ r \sin^2 \theta },
\end{equation}
and $I(r)$ joins the points $\pm \sqrt{ 1- 1/r^2 }+ i/r   $. Hence (\ref{ts5a}), (\ref{ts3}) 
 and (\ref{ts10}) yield
$$\mathfrak{m}(r, f) =
\frac1{2 \pi } \int_{- \sqrt{1-\sigma^2}}^{\sqrt{1-\sigma^2}} \log^+ | F ( t + i \sigma ) | \, dt = 
m_2(\sigma, F) + O(1) \quad \hbox{as $1/\sigma = r \to + \infty $,}
$$
since $F$ is meromorphic at
$\pm 1$.
Observe next that, still 
for $r \geq 1$ and $\sigma = 1/r$,
\begin{eqnarray*}
 \mathfrak{n}(r,f) &=&  n(\{ z \in \C :  \, |z| \geq 1 , \, |z-ir/2| \leq r/2 \}, f ) \nonumber \\
 &=& n(  \{ \zeta \in \C : \, | \zeta | \leq 1, \, {\rm Im} \, \zeta \geq \sigma  \} , F ) = n_2( \sigma, F) + O(1) . 
\end{eqnarray*}
The proof is completed by writing, with
$\lambda = 1/s$,
$$
 \mathfrak{N}(r, f) = \int_1^r \frac{  \mathfrak{n}(s, f) }{s^2} \, ds =
  \int_1^r \frac{ n_2(\lambda , F) + O(1)  }{s^2} \, ds =
 \int_\sigma^1 n_2( \lambda, F ) \, d \lambda + O(1)  .
$$
\hfill$\Box$
\vspace{.1in}

For $ \sigma \in (0, x)$ such that $F$ has no zeros or poles $\zeta = t + i \sigma$ with $-x \leq t \leq x$, write
\begin{eqnarray*}
p( \sigma ) &=& 2 \pi \,
\frac{d}{d \sigma} \left(  m_2(\sigma, F) - m_2(\sigma, 1/F)  \right)
=
\int_{-x}^x \frac{\partial}{\partial \sigma } \left( \log | F(t + i \sigma )| \right) \, dt \\
&=&
- \int_{-x}^x \frac{\partial}{\partial t } \left( \arg  F(t + i \sigma ) \right) \, dt
= -  \left( \hbox{change in $\arg F(\zeta)$ around $\partial K(\sigma)$}  \right) + O(1)\\
&=&  2 \pi (n_2(\sigma, F) - n_2(\sigma, 1/F) ) + O(1) ,
\end{eqnarray*}
using (\ref{ts10}).
Integration from $\sigma$ to $x$ yields,
as
$\sigma \to 0+$,
$$
m_2(\sigma, 1/F) - m_2(\sigma, F) =
 N_2(\sigma, F) - N_2(\sigma, 1/F) 
 + O(1),
$$
which gives (\ref{ts12}) for $a=0$, by
(\ref{101formula}),
the general case  following via standard properties of $\log^+$.

It remains to show
that
$\mathfrak{T}(r, f)$ differs from a continuous non-decreasing function by a bounded quantity
as $r \to + \infty$.
For $0 < \sigma \leq x$ set
\begin{eqnarray}
 \label{ts14}
 T_3(\sigma) &=&
 \frac1{2 \pi} \int_{-x}^x \log \sqrt{ 1+  |F(t+i \sigma )|^2}  \, dt +
N_2(\sigma, F) =
T_2(\sigma, F)
+ O(1),
 \end{eqnarray}
 by (\ref{ts10})
 and
the inequality
$0 \leq  \log \sqrt{1+y^2}
- \log^+ y \leq
\frac12 \log 2$
for $y \geq 0$
(see  \cite[p.12]{Hay2}).
Let  $\sigma \in (0, x)$ be such that $F$ has no zeros or poles $\zeta$
with 
$\zeta = t + i \sigma$ and $-x \leq t \leq x$, and hence none on $\partial K(\sigma)$.
Since  $\log \sqrt{ 1+|F|^2} $ is $C^1$  on
 $L_x$,
\begin{eqnarray}
T_3'(\sigma) &=& - n_2(\sigma,F) +
 \frac1{2 \pi} \int_{-x}^x \frac{\partial}{\partial \sigma} \left(  \log \sqrt{ 1+  |F(t+i \sigma )|^2} \right)  \, dt \nonumber \\
 &=& - n_2(\sigma, F) - 
 \frac1{2 \pi} \int_{\partial K(\sigma)}  \frac{\partial}{\partial n} \left(  \log \sqrt{ 1+  |F(t+i \sigma )|^2} \right)  \, ds
 + O(1)
 \label{ts15}
 \end{eqnarray}
 as $\sigma \to 0+$,
 in which  $\partial/\partial n$ is the outward normal derivative
 and
 $ds$ denotes integration with respect to arc length.
 Let $b_j$ be the poles of $F$ in $K(\sigma)$, with $m_j$  the corresponding multiplicities, and write
$
M(\mu)  = K(\sigma) \setminus \bigcup D(b_j, \mu),
$
for small
$\mu > 0$. Green's formula and
\cite[pp.10-11]{Hay2} yield
$$
 \int_{\partial M(\mu)}  \frac{\partial}{\partial n} \left(  \log \sqrt{ 1+  |F|^2} \right)  \, ds  =
 \int_{M(\mu)} \nabla^2  \left(  \log \sqrt{ 1+  |F|^2} \right) \, dm 
 = 2 \int_{M(\mu)} \frac{ |F'|^2 }{(1+|F|^2)^2} \, dm ,
$$
where $dm$ denotes area measure.
 Now let  $\mu \to 0+$ and observe that,  exactly as in \cite[p.11]{Hay2},
\begin{eqnarray*}
 \int_{ |\zeta-b_j| =  \mu }  \frac{\partial}{\partial n} \left(  \log \sqrt{ 1+  |F|^2} \right)  \, ds &=&
 \int_{ |\zeta-b_j| =  \mu }   \frac{m_j}{\mu} + O(1) \, ds =
 2 \pi m_j
 + o(1).
 \end{eqnarray*}
 It follows that, by (\ref{ts15}),
 as $\sigma \to 0+$,
 \begin{eqnarray*}
 \frac1\pi \int_{K(\sigma)} \frac{ |F'|^2 }{(1+|F|^2)^2} \, dm &=&
n_2(\sigma, F) +
 \frac1{2 \pi} \int_{\partial K(\sigma)}  \frac{\partial}{\partial n} \left(  \log \sqrt{ 1+  |F|^2} \right)  \, ds
= - T_3'(\sigma) + O(1)   .
\end{eqnarray*}
On combination with
(\ref{ts14}) and
Lemma
\ref{lembounded},
this delivers, as $\sigma = 1/r \to 0+$,
\begin{eqnarray*}
 \frac1\pi \int_\sigma^x \int_{K(\lambda)} \frac{ |F'|^2 }{(1+|F|^2)^2} \, dm \, d \lambda
 &=&
 T_3(\sigma) + O(1) 
 =
 \mathfrak{T}(r, f) + O(1),
\end{eqnarray*}
in which the left-hand side is evidently  non-decreasing as
$r$ increases.
 \hfill$\Box$
\vspace{.1in}


The proof of Proposition \ref{propfirsttsuji}
is complete, but it seems worth discussing briefly
the points $w \in H^+$ which contribute to $\mathfrak{m}(r, G^*) = \mathfrak{m}^*(r, G)$ in Section \ref{clunie}.
Suppose then that $w \in H^+$ and $z = |z| e^{i \phi} = s e^{i \phi}\in \Delta^+$ are related by (\ref{wdef})
and
$w$ contributes to $\mathfrak{m}(r, G^*)$, with $r \geq 1$ large. Then $w$ lies on the arc $J(r)$ in (\ref{ts1}),
and so $\zeta = -1/w $ satisfies (\ref{ts4}). Thus
$
z+1/z = 4/w 
$
has imaginary part $-4/r$, and real part lying between 
$\pm 4 \sqrt{ 1 - 1/r^2} $.
As $w$ describes the arc $J(r)$, the pre-image $z$ follows an arc in $\Delta^+$
of the polar curve
$(1/s - s) \sin \phi = 4/r$,
which passes close to $i$ when $w = ir$. If $w$ is an end-point of $J(r)$,
then $w$ lies close to $\pm 1$ and so
$z + 1/z$ is close to $\pm 4$,
from which it follows that
$z$ lies close to $\pm ( 2 - \sqrt{3} ) $.
 \hfill$\Box$
\vspace{.1in}



\section{Proof of Proposition \ref{propsecondtsuji}}\label{logderproof}

Assume now that $f$ is as in  the hypotheses of Proposition \ref{propsecondtsuji}, so that $f$ has finitely many zeros and poles in ${\rm Im} \, z \geq 0, \, |z| \leq 1$. Indeed, it may be assumed that there exists $\delta > 0$ such that $f$ has no  zeros or poles in the closed half-disc
$Y_\delta = \{ z \in \C : \, {\rm Im} \, z \geq 0, \, |z| \leq 1+\delta \} ,$
because  a rational function $S$ may be chosen such that
$f = f_1 S$ and
$f_1$ has no zeros or poles in $Y_\delta$, while
the differences
$\mathfrak{T}(r, f)- \mathfrak{T}(r, f_1)$ and
$\mathfrak{m}(r, f'/f) - \mathfrak{m}(r, f_1'/f_1)$ are bounded as $r \to + \infty$.

If the function
$g(v)$ is meromorphic on the closed disc $|v| \leq T$, with no zeros and poles on $|v| = T$, then $g'/g$ satisfies,
  for $|u| = t < T$ with $g(u) \neq 0, \infty$, 
 \begin{equation}
\label{g'gest}
\left| \frac{g'(u)}{g(u)} \right| \leq \frac{T}{\pi (T-t)^2} \int_0^{2 \pi} | \log | g( T e^{i \phi} )| | \, d \phi + 
2 \sum \frac1{|u-A_j|} ,
\end{equation} 
with the sum over all zeros and poles $A_j$ 
of $g$ in $|v| < T$, repeated according to multiplicity:
this estimate  is derived in 
 \cite{Jank} from the standard differentiated Poisson-Jensen formula \cite[p.22]{Hay2}. 


Now choose $r, s$ and $R$ with $r$ large and positive and $r < s < R \leq 2r$, such that $f$ has no zeros or poles on
the circle $C$ of centre $is/2$ and radius $s/2$.
Let $z =
r \sin \theta e^{i \theta} $ lie on the arc $J(r)$ in (\ref{ts1}), and hence inside $C$.
Applying (\ref{g'gest})  to  $g(u) = f(u+i s/2)$, with $T = s/2$, delivers
\begin{equation} 
\label{ts19}
\left| \frac{f'(z)}{f(z)} \right| \leq  \frac{s}{2\pi (s/2-|z-is/2|)^2 } 
\int_{-\pi/2}^{3 \pi/2} | \log | f( (s/2) e^{i \phi}+ is/2 )| | \, d \phi +
2 \sum \frac1{|z- B_k|} ,
\end{equation}
with the sum over all zeros and poles $B_k$ of $f$ in $|\zeta-is/2| < s/2$. 
This formula is a slightly simplified version of 
\cite[Ch.3, (3.6)]{GO}, and the subsequent steps follow closely those in \cite{GO}.

Denote by $K$  positive constants, depending  possibly on $f$ but not on $r, s$ or $R$.
Writing $\phi = 2 \eta- \pi /2$ gives
$$
(s/2) e^{i \phi}+ is/2  = \frac{is}2 ( 1 - e^{i 2 \eta} ) =  \frac{is}2 ( 2 \sin^2 \eta- 2i \sin \eta\cos \eta)
= s \sin \eta \, e^{i \eta} 
$$
and,  by (\ref{ts5a}) and the fact that  $f$ has no poles or zeros in $Y_\delta $,
\begin{eqnarray} 
I_1 &=&
  \int_{- \pi/2 }^{3 \pi/2}  \log^+ | f( (s/2) e^{i \phi}+ is/2 )|  \, d \phi
= 2 \int_{0 }^{ \pi}  \log^+ | f (s \sin \eta \, e^{i \eta}  )|  \, d \eta \nonumber \\
&\leq&  2 \int_{\arcsin(1/s) }^{ \pi - \arcsin(1/s)}  \log^+ | f (s \sin \eta \, e^{i \eta}  )|  \, d \eta+ K
\leq 4 \pi s \, \mathfrak{m}(s, f) + K .
\label{TS1}
\end{eqnarray}
The fact that $z$ lies on the arc $J(r)$ in
 (\ref{ts1}) delivers,  since $r < s < 2r$ and $1/r \leq \sin \theta \leq 1$, 
\begin{eqnarray} 
| z - is/2|^2 
&=& r (r-s) \sin^2 \theta + \frac{s^2}4  
\leq \frac{r (r-s) }{r^2} + \frac{s^2}4  < \frac{s^2}4  , \nonumber \\
s/2 - | z - is/2| &=& \frac{s^2/4 - | z - is/2|^2 }{s/2 + | z - is/2|}
\geq \frac{s-r }{r(s/2 + | z - is/2|)} 
\geq \frac{s-r}{rs}. 
\label{TS2}
\end{eqnarray} 
Proposition \ref{propfirsttsuji},  (\ref{ts19}), (\ref{TS1}) and (\ref{TS2}) now yield
\begin{eqnarray} 
\label{ts20}
\left| \frac{f'(z)}{f(z)} \right| &\leq& \frac{2 s^4 r^2}{(s-r)^2} \, (   \mathfrak{m}(s, f)  +  \mathfrak{m}(s, 1/f)  )
+ 2  \sum \frac1{|z- B_k|}   + K \nonumber \\
&\leq& \frac{4 s^4 r^2  }{(s-r)^2} \,
 \mathfrak{T}(s, f)
+ 2 \left( \sum \frac1{|z- B_k|^{1/4}} \right)^4  + K .
\end{eqnarray}
If $B_k = R_k e^{i S_k} $, with $R_k \geq 0$ and $S_k \in \R$, then
\begin{equation}
 \int_{\arcsin(1/r)}^{\pi - \arcsin(1/r)}
 \frac1{|r \sin \theta \, e^{i \theta } - B_k|^{1/4}}  \, d \theta
  \leq \frac1{r^{1/4}}
 \int_0^{2 \pi }
 \frac1{|\sin \theta \, \sin ( \theta - S_k  ) |^{1/4} }   \, d \theta
 \leq \frac{K}{ r^{1/4}} ,
 \label{I2est}
   \end{equation}
   by periodicity and
    the Cauchy-Schwarz inequality.
Next, write
\begin{equation}
\mu(r) = \int_{\arcsin(1/r)}^{\pi - \arcsin(1/r)} 
\frac{ 1  }{r \sin^2 \theta } \, d \theta = 2 \, \sqrt{1 - \frac1{r^2} }  \leq 2.
\label{muintegralest}
\end{equation}
Since $r$ is large, Jensen's inequality for an integral with respect to a probability measure yields
\begin{eqnarray}
I_2 &=& \frac1{2 \pi}  \int_{\arcsin(1/r)}^{\pi - \arcsin(1/r)}
 \log^+   \left( \sum \frac1{|r \sin \theta \, e^{i \theta } - B_k|^{1/4}} \right)  \,\frac{ d \theta }{r \sin^2 \theta }
 \nonumber
 \\
 &\leq& \frac{\mu(r)}{2 \pi} \, \int_{\arcsin(1/r)}^{\pi - \arcsin(1/r)}
 \log   \left( 1 +
 \sum \frac1{|r \sin \theta \, e^{i \theta } - B_k|^{1/4}} \right)  \,\frac{ d \theta }{\mu(r) \, r \sin^2 \theta } \nonumber \\
&\leq&  \frac{ \mu (r)}{ 2 \pi} \,
 \log \left(   \int_{\arcsin(1/r)}^{\pi - \arcsin(1/r)}
 \left( 1 +  \sum \frac1{|r \sin \theta \, e^{i \theta } - B_k|^{1/4}} \right)  \, \frac {  d \theta }{ \mu(r) \, r \sin^2 \theta } \,  \right)
 \nonumber \\
 &\leq&  \frac{1}{ \pi} \, \log
\left( 1 +  r   \sum \int_{\arcsin(1/r)}^{\pi - \arcsin(1/r)}
 \frac1{|r \sin \theta \, e^{i \theta } - B_k|^{1/4}}  \, d \theta \right)  .
 \label{I3est}
\end{eqnarray}
Combining     (\ref{ts20}),  (\ref{I2est}), (\ref{muintegralest}) and
(\ref{I3est})
delivers
\begin{eqnarray*}
\mathfrak{m}(r, f'/f) &\leq & K \left( 1 +  \log s +  \log^+ \frac1{s-r}  + \log^+ \mathfrak{T}(s, f) \right) + 4 I_2
\nonumber \\
&\leq& K \left( 1 +  \log s +  \log^+ \frac1{s-r}  + \log^+ \mathfrak{T}(s, f) + \log^+ ( \mathfrak{n}(s, f) + \mathfrak{n}(s, 1/f) ) \right).
\end{eqnarray*}
Choose  $s$  close to $(R+r)/2$: then
 (\ref{ts5a}),
 (\ref{ts7}) and
 the last estimate
 lead to
\begin{eqnarray*}
\mathfrak{T}(R, f) &\geq& \mathfrak{N}(R, f) \geq
\int_s^R \frac{\mathfrak{n}(s,f)}{t^2} \, dt = \frac{(R-s)\mathfrak{n}(s,f)}{Rs},\\
\mathfrak{m}(r, f'/f) &\leq &
 K \left( 1 +  \log R +  \log^+ \frac1{R-r}  + \log^+ \mathfrak{T}(s, f) + \log^+ \mathfrak{T}(R, f)   \right) .
 \end{eqnarray*}
 To complete the proof of
  Proposition \ref{propsecondtsuji}, take $R = 2r$ if
 $\log^+ \mathfrak{T}(t, f) = O( \log t)$ as $t \to + \infty $,
 with $R = r + ( 2 \, \mathfrak{T}(r, f))^{-1}
$ otherwise.
In the first case no exceptional set is required, while in the second case
Proposition \ref{propfirsttsuji}
gives a continuous, non-decreasing and  unbounded function $S(t)$ which differs from $\mathfrak{T}(t, f)$ by a bounded quantity.
Thus
  \cite[p.38]{Hay2}
  yields, for $r$
  outside a set of finite measure
  and $T \in \{ s, R \}$,
$$
 \mathfrak{T}(T, f) \leq S(T) + O(1)
\leq S\left( r + \frac1{S(r)} \right) + O(1)
 \leq 2S(r) + O(1) \leq 2  \mathfrak{T}( r, f)  + O(1).
 $$
 \hfill$\Box$
\vspace{.1in}

{\footnotesize

\noindent
J.K. Langley, Emeritus Professor,\\
Mathematical Sciences, University of Nottingham, NG7 2RD, UK\\
james.langley@nottingham.ac.uk

}

\end{document}